\documentclass[12pt]{article}
\usepackage[cp1251]{inputenc}
\usepackage[T2A]{fontenc}
\usepackage[english]{babel}
\usepackage{amsfonts}
\usepackage{amsmath}
\usepackage{amsmath}
\usepackage{amssymb}
\usepackage{color}

\newtheorem{theorem}{Theorem}

\newtheorem{definition}[theorem]{Definition}

\newtheorem{lemma}[theorem]{Lemma}

\begin{document}

\title{Critical Galton-Watson branching processes with countably infinitely
many types and infinite second moments \thanks{%
This work was supported from a grant of the Mathematical center in
Akademgorodok and from a grant to the Steklov International Mathematical
Center in the framework of the national project "Science" of the
Russian Federation} }
\author{V.A.Topchii\thanks{%
Mathematical center in Akademgorodok, Novosibirsk, 630090, Russia, Sobolev
Institute of Mathematics; 4 Acad. Koptyug avenue, Novosibirsk 630090,
Russia; E-mail: topchij@ofim.oscsbras.ru}, V.A.Vatutin\thanks{%
Steklov Mathematical Institute of Russian Academy of Sciences, Moscow,
Russia; 8 Gubkina St., Moscow 119991, Russia; E-mail: vatutin@mi.ras.ru},
E.E.Dyakonova\thanks{%
Steklov Mathematical Institute of Russian Academy of Sciences, Moscow,
Russia; 8 Gubkina St., Moscow 119991, Russia; E-mail: elena@mi.ras.ru} }
\date{\today }
\maketitle

\begin{abstract}
We consider an indecomposable Galton-Watson branching process with countably
infinitely many types. Assuming that the process is critical and allowing
for infinite variance of the offspring sizes of some (or all) types of
particles we describe the asymptotic behavior of the survival probability of
the process and establish a Yaglom-type conditional limit theorem for the
infinite-dimensional vector of the number of particles of all types.
\end{abstract}

\section{Definition of the process and basic properties of its mean matrix}

\label{se1}

We consider an indecomposable Galton-Watson branching process $\mathbf{Z}%
(n):= \big(Z_{j}(n)\big)_{j\in\mathbb{N}}$ with countably infinitely many
types labelled by numbers $j\in\mathbb{N}:=\{1,2,...\}$. The component $%
Z_{j}(n)$, $n\in\mathbb{N}_{0}:= \mathbb{N}\cup\{0\}$, of $\mathbf{Z}(n)$
denotes the number of type $j$ particles in the process at moment $n$. Let $%
\delta_{ij}$ be the Kronecker symbol and $\mathbf{e}_{i}:= \big(\delta_{ij}%
\big)_{j\in\mathbb{N}}$ be the vector whose $i$th component is equal to one
while the remaining are zeros. To specify the evolution of the branching
process initiated at moment $0$ by a vector $\mathbf{Z}(0)$ of individuals
of different types it is sufficient to describe the distributions of the
vectors
\[
\mathbf{Z}_{i}=\mathbf{Z}_{i}(1):=\big\{\mathbf{Z}(1)\,\big|\,\mathbf{Z}(0)=
\mathbf{e}_{i}\big\}=:\big(Z_{ij}\big)_{j\in \mathbb{N}}=\big(Z_{ij}(1)\big)%
_{j\in \mathbb{N}}.
\]

We suppose that
\begin{equation}
Z_{i}:=\sum_{j\in \mathbb{N}}Z_{ij}<\infty  \label{Finite_mu}
\end{equation}
with probability 1. Let
\[
\mathbf{Z}_{i}(n):=\big\{\mathbf{Z}(n)\,\big|\,\mathbf{Z}(0)=\mathbf{e}_{i}%
\big\}=: \big(Z_{ij}(n)\big)_{j\in \mathbb{N}}\,.
\]

Assuming that $\mathbf{s}=\big(s_{j}\big)_{j\in \mathbb{N}}\in \lbrack 0,1]^{%
\mathbb{N}}$ we specify infinite-dimensional vectors
\[
\mathbf{F}(\mathbf{s}):=\big(F_{i}(\mathbf{s})\big)_{i\in \mathbb{N}}\text{\
and \ }\mathbf{F}(n;\mathbf{s}):=\big(F_{i}(n;\mathbf{s})\big)_{i\in \mathbb{%
N}}
\]%
of the offfspring generating functions of the process with the components
\begin{eqnarray}
&&F_{i}(\mathbf{s}):=\mathbb{E}\left[ \prod_{j\in \mathbb{N}}s_{j}^{Z_{ij}}%
\right] =:\mathbb{E}\mathbf{s}^{\mathbf{Z}_{i}}=:\sum_{\mathbf{j}\in \mathbb{%
N}_{0}^{\mathbb{N}}}p_{i\mathbf{j}}\mathbf{s}^{\mathbf{j}},  \label{dfF1} \\
&&F_{i}(n;\mathbf{s}):=\mathbb{E}\left[ \prod_{j\in \mathbb{N}%
}s_{j}^{Z_{ij}(n)}\right] =:\mathbb{E}\mathbf{s}^{\mathbf{Z}_{i}(n)},
\label{dfFn}
\end{eqnarray}%
where, for any $\mathbf{j}=(j_{i})_{i\in \mathbb{N}}\in \mathbb{N}_{0}^{%
\mathbb{N}}$
\[
p_{i\mathbf{j}}:=\mathbf{P}\left( \mathbf{Z}(1)=\mathbf{j}\,\big|\,\mathbf{Z}%
(0)=\mathbf{e}_{i}\right) =\mathbf{P}\left( \mathbf{Z}_{i}(1)=\mathbf{j}%
\right) .
\]%
The probability generating functions are well defined in view of (\ref%
{Finite_mu}).

According to the branching property of the process, each type $i$ particle
belonging to the population has a unit life-length and produces at the end
of its life (independently of the prehistory of the process and the
reproduction of the others particles existing at this moment) a random
number of children specified by a vector $\mathbf{Z}_{i}$ whose distribution
is described by the generating function $F_{i}(\mathbf{s})$. This property
has the following description in terms of iterations of the offspring
generating functions:
\begin{equation}
\mathbf{F}(n+1;\mathbf{s})=\mathbf{F}(n;\mathbf{F}(\mathbf{s)})=\mathbf{F}%
\big(\mathbf{F}(n;\mathbf{s})\big)\ \mbox{for all}\ n\in \mathbb{N},
\label{df1}
\end{equation}%
where $\mathbf{F}(1;\mathbf{s}):=\mathbf{F}(\mathbf{s})$. \newline
The basic classification of the Galton-Watson branching processes with
countably infinitely many types (below we use the short abbreviation GWBP/$%
\infty $ for representatives of such processes) is given in terms of the
mean matrix
\begin{equation}
\mathbf{M}:=\big(M_{ij}\big)_{i,j\in \mathbb{N}}:=\big(\mathbb{E}Z_{ij}\big)%
_{i,j\in \mathbb{N}}=\left( \frac{\partial F_{i}(\mathbf{s})}{\partial s_{j}}%
\Big|_{\mathbf{s}=\mathbf{1}}\right) _{i,j\in \mathbb{N}}.  \label{df2}
\end{equation}

To describe such a classification in detail we recall a number of asymptotic
properties of the powers of infinite-dimensional matrices with nonnegative
elements borrowed from \cite{VJ}. We first formulate the desired properties
for an abstract matrix $M=(m_{ij})_{i,j\in \mathbb{N}}$ with $m_{ij}\geq 0$
for all $(i,j)\in \mathbb{N}^{2}$ and than use them in studying GWBP/$\infty
$ with matrix $\mathbf{M}$ specified by (\ref{df2}).

Let $M^{n}=\big(m_{ij}^{(n)}\big)_{i,j\in \mathbb{N}}$ be the $n$th power of
the infinite-dimensional matrix $M$. The matrix $M$ with nonnegative
elements is called irreducible and aperiodic if for any pair of indices $%
(i,j)$ there is an $n\in \mathbb{N}$ such that $m_{ij}^{(n)}>0$ and the
greatest common divisor of all $n\in \mathbb{N}$ such that $m_{ij}^{(n)}>0$
equals 1. According to Theorem A in \cite{VJ} for any irreducible matrix $M$
there exists a number $R\in \lbrack 0,\infty )$ such that, for any pair of
indices $(i,j)$
\begin{equation}
\lim_{n\rightarrow \infty }\left( m_{ij}^{(n)}\right) ^{1/n}=1/R.  \label{t1}
\end{equation}

The parameter $R$ is a common convergence radius of the functions
\[
\mathcal{M}_{ij}(z):=\sum_{n=0}^{\infty }m_{ij}^{(n)}z^{n},
\]
where $m_{ij}^{(0)}:=\delta _{ij},$ and $1/R$ is an analog of the maximal
(in absolute value) eigenvalue of a nonnegative matrix in the
finite-dimensional case. However, the operator specified by such a matrix is
not necessarily bounded.

It follows from (\ref{t1}) that for any pair $(i,j)$ the series
\begin{equation}
\mathcal{M}_{ij}(r)=\sum_{n=0}^{\infty }m_{ij}^{(n)}r^{n}  \label{r}
\end{equation}%
is convergent for all $0<r<R$ and is divergent for $r>R$. The case $r=R$
allows for both possibilities. Moreover, by Theorem B in \cite{VJ} the
series $\mathcal{M}_{ij}(R)$ is either convergent for all pairs $(i,j)$ or
is divergent for all $(i,j)$. In addition, the limits of all sequences $%
\left\{ m_{ij}^{(n)}R^{n},n\geq 1\right\} $ exist and ether
\begin{equation}
\lim_{n\rightarrow \infty }m_{ij}^{(n)}R^{n}=0\ \hbox{\ for all }(i,j)\in
\mathbb{N}^{2},  \label{t2}
\end{equation}%
or
\begin{equation}
\lim_{n\rightarrow \infty }m_{ij}^{(n)}R^{n}>0\ \ \hbox{\ for all }(i,j)\in
\mathbb{N}^{2}.  \label{t2b}
\end{equation}

An irreducible matrix $M$ is called

\begin{itemize}
\item ``$R$''--\textit{transient} or ``$R$''--\textit{recurrent} depending
on the convergence or divergence of the series $\mathcal{M}_{ij}(R)$ in (\ref%
{r});

\item \textquotedblleft $R$\textquotedblright --\textit{null} if (\ref{t2})
is valid and \textquotedblleft $R$\textquotedblright --\textit{positive} if (%
\ref{t2b}) holds true.
\end{itemize}

We now come back to the matrices $\mathbf{M}=\big(M_{ij}\big)_{i,j\in\mathbb{%
N}}$ and $\mathbf{M}^{n}:=\big(M_{ij}^{(n)}\big)_{i,j\in \mathbb{N}}$ of
GWBP/$\infty$, where $M_{ij}^{(n)}:=\mathbb{E}Z_{ij}(n)$, $%
M_{ij}=M_{ij}^{(1)}$.

The properties of nonnegative infinite-dimensional matrices we have listed
above allow to introduce the following natural classification of the GWBP/$%
\infty $'s (see, for instance, \cite{Sa2013}):

\begin{definition}
\label{d0} A GWBP/$\infty $ is called subcritical \{critical,
supercritical\} and transient \{recurrent, null recurrent, positively
recurrent\} in the type space, if its matrix of the mean offspring numbers $%
\mathbf{M}$ has a convergence radius $R>1$ \{$R=1,R<1$\} and is
\textquotedblleft $R$\textquotedblright -transient \{\textquotedblleft $R$%
\textquotedblright --recurrent, \textquotedblleft $R$\textquotedblright
--null recurrent, \textquotedblleft $R$\textquotedblright --positively
recurrent\}.
\end{definition}

In this paper we consider only the critical GWBP/$\infty $'s. We know by (%
\ref{t2}) or (\ref{t2b}) that if a GWBP/$\infty $ is critical then the
elements of the sequences $\left\{ M_{ij}^{(n)},n\geq 1\right\} $ either
vanishes as $n\rightarrow \infty $ for all $i$ and $j$, or have positive
limits for all $i$ and $j$. Below we analyze only the second option.

If the mean matrix $\mathbf{M}$ of a GWBP/$\infty $ is irreducible and
\textquotedblleft $1$\textquotedblright --positive then (see Theorem $D$ in
\cite{VJ}) there exist unique (up to a positive multiplier) left and right
eigenvectors $\mathbf{v}:= (v_{k})_{k\in\mathbb{N}}$ and $\mathbf{u}%
:=(u_{k})_{k\in \mathbb{N}}$ with positive components such that
\begin{equation}
\mathbf{vM}=\mathbf{v},\ \mathbf{Mu}^{T}=\mathbf{u}^{T},\ \mathbf{v}\mathbf{u%
}^{T}= \sum_{k=1}^{\infty }v_{k}u_{k}=1,\ \mathbf{v1}^{T}<\infty ,
\label{t4}
\end{equation}
and, as $n\rightarrow \infty$
\begin{equation}  \label{ed2}
M_{ij}^{(n)}\rightarrow \frac{u_{i}v_{j}}{\,\mathbf{v}\,\mathbf{u}^{T}}%
=u_{i}v_{j}
\end{equation}
for all $(i,j)\in\mathbb{N}^{2}$.

Observe that if $\mathbf{M}$ is a finite-dimensional irreducible aperiodic
matrix with its Perron root equal to 1 then $\mathbf{M}$ has, according to
the Perron-Frobenius theorem, positive left and right eigenvectors $\mathbf{v%
}$ and $\mathbf{u}$ of $\mathbf{M}$ satisfying (\ref{t4}) and, of course, $%
\mathbf{v}\mathbf{1}^{T}<\infty $ in this case. This estimate is, in
general, not valid for irreducible and ``$1$''--positive
infinite-dimensional matrices. We require in the paper finiteness of the
scalar product $\mathbf{v}\mathbf{1}^{T}$ (matrices with $\mathbf{v}\mathbf{1%
}^{T}<\infty $ are called irreducible with finite iterate coefficients, see,
for instance, \cite{R2014}) leaving the case $\mathbf{v}\mathbf{1}%
^{T}=\infty $ for the future investigations.

We now introduce an important definition which incorporates major
re\-strictions on the properties of the matrix $\mathbf{M}$.

Set $M_{i}:=\mathbb{E}Z_{i}=\sum_{j\in \mathbb{N}}M_{ij}$.

\begin{definition}
\label{d2} We say that the mean matrix $\mathbf{M}=\left(\mathbb{E}%
Z_{ij}\right)_{i,j\in\mathbb{N}} $ of a critical GWBP/$\infty$ belongs to a
class $\mathcal{M}_{1}$, and write $\mathbf{M}\in\mathcal{M}_{1}$, if

\begin{itemize}
\item[$\mathrm{(i)}$] $\mathbf{M}$ is irreducible and aperiodic,
\textquotedblleft $1$\textquotedblright --positive and \textquotedblleft $1$%
\textquotedblright --recurrent;

\item[$\mathrm{(ii)}$] the left eigenvector $\mathbf{v}$ of $\mathbf{M}$ has
 $L_{1}$-norm equal to 1: $\Vert \mathbf{v}\Vert _{1}=\mathbf{v}\mathbf{1}%
^{T}=1$, and the right eigenvector $\mathbf{u}$ has finite $L_{\infty
} $-norm: $U:=\Vert \mathbf{u}\Vert _{\infty }=\sup_{i\in \mathbb{N}%
}u_{i}<\infty $;

\item[$\mathrm{(iii)}$] $\lim\limits_{N\rightarrow \infty }\sup\limits_{i\in
\mathbb{N}}M_{i}^{-1}\sum\limits_{j>N}M_{ij}=0$ and $\lim\limits_{K%
\rightarrow \infty }\sup\limits_{i\in \mathbb{N}}M_{i}^{-1}\mathbb{E}\left[
Z_{i}\,;\,Z_{i}>K\right] =0$.
\end{itemize}

We say that $\mathbf{M}$ belongs to a subclass $\mathcal{M}_{1}^{0}\subset
\mathcal{M}_{1}$ if, additionally,

\begin{itemize}
\item[$\mathrm{(iv)}$] there exist $m\in\mathbb{N}$ and $c,C\in\mathbb{R}%
^{+} $ such that
\begin{equation}  \label{ed03}
M_{ij}<Cu_{i}v_{j},\quad M_{1j}^{(m)}>cv_{j},\ \forall \,i,j\in \mathbb{N}.
\end{equation}
\end{itemize}
\end{definition}

We suppose, without loss of generality, that if $\mathbf{M}\in \mathcal{M}%
_{1}$ then $\mathbf{v1}^{T}=1$ and, therefore, relations (\ref{t4}) and
Condition $\mathrm{(ii)}$ have the component-wise representations
\begin{equation}
\begin{array}{l}
\sum\limits_{j\in \mathbb{N}}M_{ij}u_{j}=u_{i};\ \sum\limits_{j\in \mathbb{N}%
}v_{j}M_{ji}=v_{i};\ \sum\limits_{j\in \mathbb{N}}v_{j}u_{j}=1; \\
\sum\limits_{j\in \mathbb{N}}v_{j}=1;\ \sup\limits_{i\in \mathbb{N}%
}u_{i}<\infty ;\ v_{j}u_{j}>0,\ \forall j\in \mathbb{N}.%
\end{array}
\label{ed3}
\end{equation}

Observe that Condition $\mathrm{(iii)}$ has rather transparent meaning. Its
first part extracts from all critical GWBP/$\infty $'s those processes in
which particles of all types produce with a high probability particles of
types with relatively small labels. Thus, our model is, in a sense, close
to the so-called lower Hessenberg branching processes \cite{BH2018} in which particles
of type $i$ may produce particles of types $j\leq i+1$ only. The second part
of Condition $\mathrm{(iii)}$ prevents existence of very productive
particles.

We need a number of auxiliary functions related to the generating functions $%
F_{i} (n;\mathbf{s})$, $i\in \mathbb{N}$, of the GWBP/$\infty$ $\mathbf{Z}%
_{i}(n)$:
\begin{eqnarray*}
&&F_{ij}(s):=\mathbb{E}s^{Z_{ij}};\ \ Q_{i}(n;\mathbf{s}):=1-F_{i}(n;\mathbf{%
s})=1- \mathbb{E}\mathbf{s}^{\mathbf{Z}_{i}(n)}; \\
&&\mathbf{Q}(n;\mathbf{s}):=\big(Q_{i}(n;\mathbf{s})\big)_{i\in \mathbb{N}}=
\mathbf{1}-\mathbf{F}(n;\mathbf{s});\ \ \mathbf{Q}(\mathbf{s}):=\mathbf{Q}%
(1; \mathbf{s})=\mathbf{1}-\mathbf{F}(\mathbf{s}); \\
&&Q_{i}(n):=Q_{i}(n;\mathbf{0})=\mathbb{P}\big(\boldsymbol{Z}(n)\neq 0\big|
\boldsymbol{Z}(0)=\mathbf{e}_{i}\big); \\
&&\mathbf{Q}(n):=\mathbf{Q}(n;\mathbf{0}),\ \ q(n;\mathbf{s}):=\mathbf{vQ}%
^{T}(n; \mathbf{s});\ \ q(n):=q(n;\mathbf{0}).
\end{eqnarray*}

For $x\geq 0$ and $U=\sup_{i\in \mathbb{N}}u_{i}$ introduce the function
\begin{equation}
\Phi (x):=\left\{
\begin{array}{lll}
x-\mathbf{vQ}^{T}(\mathbf{1}-x\mathbf{u}) & \text{if} & 0\leq xU\leq 1, \\
x-\mathbf{vQ}^{T}(\mathbf{0}) & \text{if} & xU>1.%
\end{array}%
\right.  \label{V2}
\end{equation}

We now may formulate the main result of the paper.

\begin{theorem}
\label{t02} Let $\left\{ \mathbf{Z}_{i}(n),i\in \mathbb{N}\right\} $ be a
critical GWBP/$\infty$ with mean matrix $\mathbf{M}\in\mathcal{M}_{1}^{0}$
and
\begin{equation}  \label{vat1}
\Phi (x)=x^{\alpha+1}\ell (x),
\end{equation}
where $\alpha\in(0,1]$ and $\ell(x)$ is a slowly varying function as $%
x\rightarrow +0$.

Then

1) for some slowly varying function $\ell_{1}(n)$
\begin{equation}  \label{ted22}
q(n)=n^{-1/\alpha }\ell _{1}(n)
\end{equation}
as $n\rightarrow \infty $;

2) for any $i\in\mathbb{N}$
\begin{equation}  \label{ted2222}
Q_{i}(n)=\mathbb{P}\big(\mathbf{Z}(n)\neq\mathbf{0}\big|\mathbf{Z}(0)=
\mathbf{e}_{i}\big)\sim u_{i}n^{-1/\alpha}\ell_{1}(n)
\end{equation}
as $n\rightarrow \infty $;

3) for each vector $\boldsymbol{\lambda }=\left( \lambda _{k}\right) _{k\in
N}$ with bounded coordinates and each $i\in\mathbb{N}$
\begin{equation}  \label{tv1_new}
\lim_{n\rightarrow \infty }\mathbb{E}\left[e^{-(\boldsymbol{\lambda },%
\mathbf{Z} (n))q(n)}\big|\mathbf{Z}(n)\neq \mathbf{0},\mathbf{Z}(0){=}%
\mathbf{e}_{i} \right] =1-\big(1+(\mathbf{v},\boldsymbol{\lambda })^{-\alpha}%
\big)^{-1/\alpha}.
\end{equation}
In particular, for each vector $(z_{1},\ldots,z_{m})\in\mathbb{R}_{+}^{m}$
the limit
\begin{equation}  \label{ted11}
G_{m}(z_{1},\ldots,z_{m}){:=}\lim_{n\rightarrow \infty }\mathbb{P}\left(
Z_{j}(n)q(n){\leq} z_{j},j{=}1,{...},m\big|\mathbf{Z}(n){\neq}\mathbf{0},%
\mathbf{Z}(0)=\mathbf{e}_{i}\right)
\end{equation}
exists and is independent of $i$.
\end{theorem}

\textbf{Remark}. It follows from (\ref{tv1_new}) that, as $%
n\rightarrow\infty $
\[
\left\{ Z(n)q(n)\big|\mathbf{Z}(n)\neq\mathbf{0},\mathbf{Z}(0)=\mathbf{e}%
_{i}\right\} \overset{d}{\rightarrow}\xi\mathbf{v},
\]
where
\[
\mathbb{E}e^{-t\xi }=1-\big(1+t^{-\alpha}\big)^{-1/\alpha},\ t\geq 0.
\]

We note that Kolmogorov \cite{Ko1938} was the first who investigated the
asymptotic behavior of a single-type critical Galton-Watson process. His
work was followed by the celebrate Yaglom article \cite{Ya1947} who studied
the distribution of the number of particles in a single-type critical
Galton-Watson process given its survival for a long time. Joffe and Spitzer
\cite{JH67} extended these results to the case of multi-type critical
indecomposable Galton-Watson processes. All these papers required finiteness
of the second moments of the reproduction laws of the number of particles.

Zolotarev \cite{Zol57}, assuming that the variance for the offspring
reproduction law of particles may be infinite, had found an asymptotic
representation for the survival probability of a single-type continuous-time
critical branching processes and proved a Yaglom-type theorem for this case.
Zolotarev's results were complemented by Slack \cite{Sl68}, \cite{Sl72} who
generalized Kolmogorov's and Yaglom's theorems to the case when the
offspring generating function of a critical Galton-Watson process has the
form
\[
F(s)=s+\left( 1-s\right) ^{1+\alpha }\ell (1-s),
\]%
where $\alpha \in (1,2]$ and $\ell (x)$ is a slowly varying function as $%
x\rightarrow +0$.

Slack's theorems were independently and almost simultaneously extended to
the multi-type indecomposable setting by Vatutin \cite{V1977} and Goldstein
and Hoppe \cite{GH}. The main assumption of these two papers is just our
condition~(\ref{vat1}) formulated in terms of the critical Galton-Watson
processes with finite number of types. Thus, Theorem \ref{t02} is a natural
generalization of the main results of \cite{V1977} and \cite{GH} to the GWBP/%
$\infty $'s.

There are several published results for GWBP/$\infty $ (see, for example,
\cite{AKan}, \cite{BDH2017}, \cite{BH2018}, \cite{HLN}, \cite{Kes89} and
\cite{Sa2013}). Sagitov's article \cite{Sa2013} is the most relevant to our
paper. The author analyzed there the case of linear-fractional offspring
generating functions. He has established, along with other results, an
asymptotic representation for the survival probability of a critical GWBP/$%
\infty $ and proved a Yaglom-type conditional limit theorem for such
processes. Theorem \ref{t02} of our paper extends the mentioned Sagitov
result in two directions. First, we consider the general form of the
reproduction generating functions of particles and, second, we do not assume
finiteness of the second moments for the offspring numbers of particles.

The paper is organized as follows. In Section \ref{se2} we prove a number of
statements describing properties of the offspring generation functions of
the GWBP/$\infty $'s and show that the dichotomy property, which states that
with probability 1 the population either becomes extinct or drifts to
infinity, holds for the processes meeting the conditions of Theorem \ref{t02}.

One of the basic assumptions of Theorem \ref{t02} is condition (\ref{vat1})
expressed in terms of the eigenvectors $\mathbf{v}$ and $\mathbf{u}$ of the
mean matrix $M$ and a single variable $x$. The goal of Section \ref{SecRatio}
is to demonstrate that properties of iterations of the offspring generating
functions depending on the unbounded number of arguments may be reduced to
considering some function which depends on a single argument only. To this
aim we prove a Ratio Theorem \ref{t01_splitted} showing that the functions $%
Q_{i}(n;\mathbf{s})$ may be well approximated by $u_{i}q(n;\mathbf{s})$ for
all $i\in \mathbb{N}$. This approximation allows us to complete the proof of
Theorem \ref{t02} by the methods similar to those used in \cite{V1977} for
the case of Markov branching processes with finite number of types.

\section{Properties of generating functions}

\label{se2}

We prove in this section a number of statements describing properties of the
offspring generating functions of a critical GWBP/$\infty$. Some of these
statements look evident for the Galton-Watson processes with finite number
of types. However, certain efforts and restrictions are needed to check
their validity for the infinite type case. The first result of such a kind
is the following lemma.

\begin{lemma}
\label{loc10} If $\mathbf{M}\in\mathcal{M}_{1}$ then
\begin{equation}  \label{DefMatrix0}
\liminf_{n\rightarrow\infty}\mathbf{M}^{n}\mathbf{1}^{T}\in\mathbb{R}^{%
\mathbb{N}}.
\end{equation}

If $\mathbf{M\in }\mathcal{M}_{1}^{0}$ then there exists a constant $C\in
(0,\infty )$ such that
\begin{equation}
M_{i}^{(n)}:=\sum_{j\in \mathbb{N}}M_{ij}^{(n)}\leq Cu_{i}\leq CU=:\mathfrak{%
m}<\infty  \label{ed04}
\end{equation}%
for all $i$ and $n$ belonging to the set $\mathbb{N}$.
\end{lemma}

\textbf{Remark.} Observe the difference between the estimates (\ref%
{DefMatrix0}) and (\ref{ed04}). For the first case the liminf of the
row-wise sums of elements is finite while for the second one the sums are
uniformly bounded. Clearly, the second statement is not a consequence of the
first one.

\textbf{Proof of Lemma \ref{loc10}.} Fix an $\varepsilon \in (0,0.5)$ and,
using Condition $\mathrm{(ii)}$ of Definition \ref{d2} select a positive
integer $N=N(\varepsilon )$ such that
\begin{equation}
\sup_{k\in \mathbb{N}}\sum_{j\geq N}M_{kj}\leq \varepsilon .  \label{le2}
\end{equation}

Fix now an $i\in \mathbb{N}$. Recalling the conditions $\mathbf{v1}^{T}=1$, $%
\mathbf{u}>\mathbf{0}$, $\mathbf{v}>\mathbf{0}$, and the limiting relation (%
\ref{ed2}), we conclude that for $\delta :=\sum_{j\geq N}v_{j}$ there exists
$n_{0}=n_{0}(i,N)$ such that the estimate
\begin{equation}
\sum_{j<N}M_{ij}^{(n)}\leq (1+\delta )u_{i}\sum_{j<N}v_{j}=(1-\delta
^{2})u_{i}\leq u_{i}  \label{le21}
\end{equation}%
is valid for all $n\geq n_{0}$.

Using (\ref{le2}) and (\ref{le21}) for $n\geq n_{0}$ gives
\[
M_{i}^{(n)}=\sum_{j<N}M_{ij}^{(n)}+\sum_{j\geq N}\sum_{k\in \mathbb{N}%
}M_{ik}^{(n-1)} M_{kj}\leq u_{i}+\varepsilon M_{i}^{(n-1)}
\]
or, for $n\geq 1$
\begin{equation}  \label{le31}
M_{i}^{(n+n_{0})}\leq u_{i}\sum_{l=0}^{n-1}\varepsilon ^{l}+\varepsilon
^{n}M_{i}^{(n_{0})}.
\end{equation}

Note that $M_{i}=\mathbb{E}Z_{i}\leq\mathbb{E}\left[Z_{i}\,;\,Z_{i}>K\right]+K$. So the second part of Condition $\mathrm{(iii)}$ provides existence of a
constant $W<\infty$ such that
\begin{equation}
\sup_{k\in \mathbb{N}}M_{i}\leq W .  \label{rA}
\end{equation}

Hence we deduce the following estimate which is valid for all $i\in\mathbb{N}
$:
\[
M_{i}^{(n_{0})}=\sum_{j\in \mathbb{N}}\sum_{k\in\mathbb{N}%
}M_{ik}^{(n_{0}-1)}M_{kj} \leq WM_{i}^{(n_{0}-1)}\leq W^{n_{0}}.
\]

This fact combined with (\ref{le31}) completes the proof of (\ref{DefMatrix0}%
).

To check the validity of the second statement of Lemma \ref{loc10} observe
that $\mathbf{vM}^{n}=\mathbf{v}$, $\mathbf{M}^{n}\mathbf{u}^{T}=\mathbf{u}%
^{T}$ for all $n\in \mathbb{N}$. Hence, using (\ref{ed03}) we conclude that,
for all $i,n\in \mathbb{N}$
\begin{equation}
M_{i}^{(n)}=\sum_{j\in \mathbb{N}}\sum_{k\in \mathbb{N}}M_{ik}^{(n-1)}M_{kj}%
\leq C\sum_{j\in \mathbb{N}}\sum_{k\in \mathbb{N}%
}M_{ik}^{(n-1)}u_{k}v_{j}=Cu_{i}.  \label{ed040}
\end{equation}

The last implies (\ref{ed04}), since $\left\Vert \mathbf{u}%
\right\Vert_{\infty }=U<\infty$ by Condition $\mathrm{(ii)}$.

Lemma \ref{loc10} is proved.

\begin{lemma}
\label{L_Dichotom0} If $\mathbf{M}\in \mathcal{M}_{1}$ and $\mathbf{F}(%
\mathbf{s})\neq \mathbf{Ms}$ then for each $i\in \mathbb{N}$ there exists $%
n=n(i)$ such that
\begin{equation}
F_{i}(n;\mathbf{0})=\mathbb{P}\left( \big\|\mathbf{Z}_{i}(n)\big\|%
_{1}=0\right) >0.  \label{Posit}
\end{equation}
\end{lemma}

\textbf{Proof}. Assume the contrary that there exists $i\in \mathbb{N}$ such
that
\begin{equation}
F_{i}(n;\mathbf{0})=0\mbox{ for all }n\in \mathbb{N}.  \label{ZeroValue}
\end{equation}

We split the set of types of the process into two parts $\mathcal{T}_{1}$
and $\mathcal{T}_{2}$. We assign type $i$ to the class $\mathcal{T}_{1}$ if
condition (\ref{ZeroValue}) holds and to the class $\mathcal{T}_{2}$ if
relation (\ref{Posit}) is valid for some $n=n(i)$.

Denote
\[
\widehat{M}_{i}:=\sum_{k\in \mathcal{T}_{1}}\mathbb{E}Z_{ik}=\sum_{k\in
\mathcal{T}_{1}}M_{ik}.
\]

Observe that $\widehat{M}_{i}\geq 1$ for all $i\in \mathcal{T}_{1}$. By
induction it is easy to show that, for any $i\in \mathcal{T}_{1}$
\begin{eqnarray*}
\widehat{M}_{i}^{(n+1)}&:=&\sum_{k\in \mathcal{T}_{1}}\mathbb{E}%
Z_{ik}(n+1)=\sum_{k\in \mathcal{T}_{1}}M_{ik}^{(n+1)} \\
&=&\sum_{k\in \mathcal{T}_{1}}\sum_{r\in \mathbb{N}}M_{ir}^{(n)}M_{rk}=%
\sum_{r\in \mathbb{N}}M_{ir}^{(n)}\widehat{M}_{r} \\
&\geq &\sum_{r\in \mathcal{T}_{1}}M_{ir}^{(n)}\widehat{M}_{r}\geq \sum_{r\in
\mathcal{T}_{1}}M_{ir}^{(n)}=\widehat{M}_{i}^{(n)}\geq 1.
\end{eqnarray*}

To go further we need to separately consider the cases $\mathcal{T}_{2}\neq
\emptyset$ and $\mathcal{T}_{2}=\emptyset$.

1) Assume first that $\mathcal{T}_{2}\neq \emptyset $. Since $\mathbf{M}$ is
an irreducible matrix, it follows that one can find $i_{0}\in \mathcal{T}%
_{1} $ and $j_{0}\in \mathcal{T}_{2}$ such that $M_{i_{0}j_{0}}=:\Delta >0$
and, therefore,
\[
M_{i_{0}}=\widehat{M}_{i_{0}}+\sum_{k\in\mathcal{T}_{2}}M_{i_{0}k}\geq
1+\Delta.
\]

For the same reason there exists $n_{0}$ such that
\[
M_{j_{0}i_{0}}^{(n_{0})}=\mathbb{E}Z_{j_{0}i_{0}}(n_{0})\geq \mathbb{P}%
\left( Z_{j_{0}i_{0}}(n_{0})>0\right) >0.
\]%
Setting $\Delta _{1}:=M_{i_{0}j_{0}}M_{j_{0}i_{0}}^{(n_{0})}>0$ we have
\begin{eqnarray*}
\widehat{M}_{i_{0}}^{(n_{0}+1)} &=&\sum_{k\in \mathcal{T}%
_{1}}M_{i_{0}k}^{(n_{0}+1)}=\sum_{k\in \mathcal{T}_{1}}\sum_{r\in \mathbb{N}%
}M_{i_{0}r}M_{rk}^{(n_{0})} \\
&=&\sum_{r\in \mathbb{N}}M_{i_{0}r}\widehat{M}_{r}^{(n_{0})}=\sum_{r\in
\mathcal{T}_{1}}M_{i_{0}r}\widehat{M}_{r}^{(n_{0})}+\sum_{r\in \mathbb{N}%
\backslash \mathcal{T}_{1}}M_{i_{0}r}\widehat{M}_{r}^{(n_{0})} \\
&\geq &\sum_{r\in \mathcal{T}%
_{1}}M_{i_{0}r}+M_{i_{0}j_{0}}M_{j_{0}i_{0}}^{(n_{0})}=\widehat{M}%
_{i_{0}}+\Delta _{1}\geq 1+\Delta _{1}.
\end{eqnarray*}

By the same arguments we conclude that, for any $q\in \mathbb{N}$
\begin{eqnarray*}
\widehat{M}_{i_{0}}^{(\left( q+1\right) n_{0}+q+1)} &=&\sum_{k\in \mathcal{T}%
_{1}}\sum_{r\in \mathbb{N}}M_{i_{0}r}^{(qn_{0}+q)}M_{rk}^{(n_{0}+1)} \\
&\geq &\sum_{r\in \mathcal{T}_{1}}M_{i_{0}r}^{(qn_{0}+q)}\widehat{M}%
_{r}^{(n_{0}+1)} \\
&\geq &\widehat{M}_{i_{0}}^{(qn_{0}+q)}+\Delta
_{1}M_{i_{0}i_{0}}^{(qn_{0}+q)}\geq \Delta
_{1}\sum_{t=1}^{q}M_{i_{0}i_{0}}^{(tn_{0}+t)}.
\end{eqnarray*}

Since $\mathbf{M}$ is a \textquotedblleft $1$\textquotedblright -recurrent
and \textquotedblleft $1$\textquotedblright -positive matrix, we have by (\ref{ed2})
\[
\sum_{t=1}^{\infty }M_{i_{0}i_{0}}^{(tn_{0}+t)}=\infty .
\]%
Hence it follows that, as $q\rightarrow \infty $
\[
\widehat{M}_{i_{0}}^{(qn_{0}+q)}\rightarrow \infty
\]%
contradicting (\ref{ed04}). Thus, if (\ref{ZeroValue}) holds then $\mathcal{T%
}_{2}$ may be only an empty set.

2) Assume now that (\ref{ZeroValue}) holds and $\mathcal{T}_{2}=\emptyset $.
In this case $F_{i}(1;\mathbf{0})=F_{i}(\mathbf{0})=0$ for all $i\in \mathbb{N}$. Therefore, the offspring generating functions may be written for all $i\in
\mathbb{N}$ as
\[
F_{i}(\mathbf{s})=\sum_{\boldsymbol{j}\in \mathbb{N}_{0}^{\mathbb{N}%
}\backslash \{\mathbf{0}\}}p_{i\boldsymbol{j}}\mathbf{s}^{\boldsymbol{j}%
},\quad F_{i}(\mathbf{1})=\sum_{\boldsymbol{j}\in \mathbb{N}_{0}^{\mathbb{N}%
}\backslash \{\mathbf{0}\}}p_{i\boldsymbol{j}}=1.
\]%
Hence,
\[
M_{i}=\widehat{M}_{i}=\sum_{\boldsymbol{j}\in \mathbb{N}_{0}^{\mathbb{N}%
}\backslash \left\{ \mathbf{0}\right\} }\left\Vert \mathbf{j}\right\Vert
_{1}p_{i\boldsymbol{j}}\geq 1
\]%
for all $i\in \mathbb{N}$, where $\left\Vert \mathbf{j}\right\Vert
_{1}=j_{1}+j_{2}+\ldots$ . It is not difficult to see that the case $M_{i}=1$
for all $i\in \mathbb{N}$ is possible only if $p_{i\boldsymbol{j}}=0$ for
all $\left\Vert \mathbf{j}\right\Vert _{1}\geq 2$, i.e. for $\mathbf{F}(%
\mathbf{s})\equiv \mathbf{Ms}$, which is not allowed by our assumptions.
Consequently, there exist $i_{0}$ and $\mathbf{j}_{0}$ with $\left\Vert
\mathbf{j}_{0}\right\Vert _{1}\geq 2$ such that $p_{i_{0}\boldsymbol{j}%
_{0}}>0$. Clearly, $M_{i_{0}}>1$ in this case.

Further, there exists $n_{0}$ such that $M_{i_{0}i_{0}}^{(n_{0})}=\Delta
_{1}>0$. Repeating now almost literally the arguments used to analyze the
case $\mathcal{T}_{2}\neq \emptyset $ we conclude that $M_{i_{0}}^{(n)}%
\rightarrow \infty $ as $n\rightarrow \infty $. This contradicts to the
uniform boundness of $M_{i}^{(n)}$ for all $i\in \mathbb{N}$. Thus, the case
$\mathcal{T}_{2}=\emptyset $ is also impossible under the assumption (\ref%
{ZeroValue}). The obtained contradiction proves (\ref{Posit}).

Lemma \ref{L_Dichotom0} is proved.

The next lemma is a refinement of Lemma \ref{L_Dichotom0}.

\begin{lemma}
\label{L_Dichotom} If $\mathbf{M}\in\mathcal{M}_{1}$ and $\mathbf{F}(\mathbf{%
s})\neq \mathbf{Ms}$ then
\[
\mathbb{P}\left( \lim_{n\rightarrow \infty }\big\|\mathbf{Z}_{i}(n)\big\|%
_{1}=0\right) =1
\]
for all $i\in\mathbb{N}$.
\end{lemma}

\textbf{Proof}. First observe that
\begin{equation}
\mathbb{P}\left( \lim_{n\rightarrow \infty }\big\|\boldsymbol{Z}_{i}(n)\big\|%
_{1}=\infty \right) =0  \label{Zero}
\end{equation}%
for all $i\in \mathbb{N}$. Indeed, if it would be not the case then
\[
\limsup_{n\rightarrow \infty }M_{i}^{(n)}=\limsup_{n\rightarrow \infty }%
\mathbb{E}\big\|\boldsymbol{Z}_{i}(n)\big\|_{1}=\infty
\]%
for some $i$, contradicting (\ref{ed040}).

Thus, to prove the lemma it is sufficient to establish that under our
conditions the process obeys the so-called dichotomy property (see, for
instance, \cite{Ha1966}):
\[
\mathbb{P}\left( \lim_{n\rightarrow \infty }\big\|\boldsymbol{Z}_{i}(n)\big\|%
_{1}=\infty \right) +\mathbb{P}\left( \lim_{n\rightarrow \infty }\big\|%
\boldsymbol{Z}_{i}(n)\big\|_{1}=0\right) =1.
\]

It is shown in \cite{Tez2002} (see Condition 2.1 and the proof of
Proposition 2.2 there) that if, for all $k\in \mathbb{N}$, there exist an
index $m_{k}$ and a positive number $d_{k}$ such that
\begin{equation}  \label{dich}
\inf_{i\in \mathbb{N}}\mathbb{P}\left( \|\mathbf{Z}_{i}(m_{k})\|_{1}=0\big|%
1\leq \|\mathbf{Z}_{i}(1)\|_{1}\leq k\right) \geq d_{k}
\end{equation}
then the respective process possesses the dichotomy property.

Let us check that (\ref{dich}) is valid if there exist an index $m_{0}$ and
a real number $d_{0}\in (0,1)$, such that
\begin{equation}  \label{dich1}
\inf_{i\in \mathbb{N}}\mathbb{P}\left( \big\|\mathbf{Z}_{i}(m_{0})\big\|%
_{1}=0\right) \geq d_{0}.
\end{equation}

Indeed, take $\mathbf{r}=(r_{j})_{j\in \mathbb{N}}\in \mathbb{N}^{\mathbb{N}%
}_{0}$ and introduce the set of events $\mathcal{A}_{i,m}:=\{\|\mathbf{Z}%
_{i}(m)\|_{1}=0\}$ and
\[
\mathcal{B}_{i,k}:=\{1\leq \|\mathbf{Z}_{i}(1)\|_{1}\leq k\} =\sum_{1\leq
\Vert \mathbf{r}\Vert _{1}\leq k}\{Z_{ij}=r_{j}\}_{j\in \mathbb{N}%
}=:\sum_{1\leq \Vert \mathbf{r}\Vert _{1}\leq k}\mathcal{B}_{i\mathbf{r}}.
\]

Since
\[
\mathbb{P}(\mathcal{A}_{i,m_{0}+1}|\mathcal{B}_{i\mathbf{r}})=\prod_{j\in
\mathbb{N}}\mathbb{P}^{r_{j}}\big(\mathcal{A}_{j,m_{0}}\big)\geq d_{0}^{k}
\]
for all $i\in \mathbb{N}$, it follows by the total probability formula that
\begin{eqnarray*}
\mathbb{P}\big(\mathcal{A}_{i,m_{0}+1}\big|\mathcal{B}_{i,k}\big)&=&\frac{ %
\textstyle{\sum\limits_{1\leq\Vert\mathbf{r}\Vert _{1}\leq k}\mathbb{P}(
\mathcal{A}_{i,m_{0}+1}\mathcal{B}_{i\mathbf{r}})}}{\textstyle{\mathbb{P}(
\mathcal{B}_{i,k})}} \\
&=&\frac{\textstyle{\sum\limits_{1\leq \Vert \mathbf{r}\Vert _{1}\leq k}
\mathbb{P}(\mathcal{A}_{i,m_{0}+1}|\mathcal{B}_{i\mathbf{r}})\mathbb{P}(
\mathcal{B}_{i,\mathbf{r}})}}{\textstyle{\mathbb{P}(\mathcal{B}_{i,k})}}\geq
d_{0}^{k}.
\end{eqnarray*}
This proves (\ref{dich}) with $m_{k}=m_{0}+1$ and $d_{k}=d_{0}^{k}$.

We now show that the estimate (\ref{dich1}) indeed holds under the
conditions of Lemma \ref{L_Dichotom} .

According to the first part of Condition $\mathrm{(iii)}$, for each $%
\varepsilon \in (0,1)$ there exists $N=N(\varepsilon )$ such that
\begin{equation}
\sup_{i\in \mathbb{N}}\mathbb{P}\left( \sum_{j>N}Z_{ij}>0\right) \leq
\sup_{i\in \mathbb{N}}\sum\limits_{j>N}M_{ij}\leq \varepsilon \sup_{i\in
\mathbb{N}}M_{i}\leq \varepsilon \mathfrak{m}.  \label{dich2}
\end{equation}

We split types of particles into two groups $\mathrm{T}_{1}:=\{j\leq N\}$
and $\mathrm{T}_{2}:=\{j>N\}$ and consider the sets
\[
\mathcal{A}_{i,m}^{\mathrm{T}_{1}}:=\left\{ \sum_{j\leq
N}Z_{ij}(m)=0\right\} \text{ \ \ and }\quad \mathcal{A}_{i,m}^{\mathrm{T}%
_{2}}:=\left\{ \sum_{j>N}Z_{ij}(m)=0\right\} .
\]

By (\ref{dich2}) we have
\begin{equation}
P_{2}(1):=\inf_{i\in \mathbb{N}}\mathbb{P}\left( \mathcal{A}_{i,1}^{\mathrm{T%
}_{2}}\right) \geq 1-\varepsilon \mathfrak{m}.  \label{dich3}
\end{equation}

In view of the second part of Condition $\mathrm{(iii)}$, for each $%
\varepsilon \in (0,1)$ there exists $K=K(\varepsilon )$ such that
\begin{equation}
\sup\limits_{i\in \mathbb{N}}\mathbb{E}\left[ Z_{i}\,;\,Z_{i}>K\right] \leq
\varepsilon M_{i}\leq \varepsilon \mathfrak{m}.  \label{dich4}
\end{equation}

Similarly to (\ref{dich3}) we have
\begin{equation}
\inf\limits_{i\in \mathbb{N}}\mathbb{P}(Z_{i}\leq K)\geq 1-\varepsilon
\mathfrak{m}.  \label{dich5}
\end{equation}

By Lemma \ref{L_Dichotom0} for each fixed $i$ there exists $n(i)$ such that $%
Q_{i}(n(i);\mathbf{0})<1$. Thus, there exist $n_{0}$ and $\theta \in (0,1)$
such that, for all $i\in \mathrm{T}_{1}$ and $n\geq n_{0}$
\begin{equation}  \label{dich6}
Q_{i}(n;\mathbf{0})\leq Q_{i}(n_{0};\mathbf{0})\leq 1-\theta <1
\end{equation}
or, for all $n\geq n_{0}$
\begin{equation}  \label{dich7}
P_{1}(n):=\inf_{i\in \mathrm{T}_{1}}\mathbb{P}\big(\mathcal{A}_{i,n}\big)%
>\theta >0.
\end{equation}

Note that if $\mathcal{B}_{1}$ and $\mathcal{B}_{2}$ are two events such
that $\mathbb{P}(\mathcal{B}_{1})>1-\sigma _{1}$ and $\mathbb{P}(\mathcal{B}%
_{2})>1-\sigma _{2}$ for some constants $\sigma _{1},\sigma _{2}\in (0,1)$,
then $\mathbb{P}(\mathcal{B}_{1}\mathcal{B}_{2})>1-\sigma _{1}-\sigma _{2}$.
Using this simple observation and recalling (\ref{dich3}) and (\ref{dich5})
gives
\begin{eqnarray*}
\inf_{i\in \mathbb{N}}\mathbb{P}\big(\mathcal{A}_{i,n_{0}+1}\big) &\geq
&\inf_{i\in \mathbb{N}}\mathbb{P}\big(\mathcal{A}_{i,n_{0}+1};Z_{i}\leq K,%
\mathcal{A}_{i,1}^{\mathrm{T}_{2}}\big) \\
&\geq &\inf_{i\in \mathbb{N}}\mathbb{P}\big(\mathcal{A}_{i,n_{0}+1}\big|%
Z_{i}\leq K,\mathcal{A}_{i,1}^{\mathrm{T}_{2}}\big)(1-2\varepsilon \mathfrak{%
m}) \\
&\geq &\inf_{i\in \mathrm{T}_{1}}\mathbb{P}^{K}\big(\mathcal{A}_{i,n_{0}}%
\big)(1-2\varepsilon \mathfrak{m})\geq \theta ^{K}(1-2\varepsilon \mathfrak{m%
}).
\end{eqnarray*}%
Selecting $\varepsilon \in \big(0,0.5\mathfrak{m}^{-1}\big)$ we justify (\ref%
{dich1}) and complete the proof of Lemma~\ref{L_Dichotom}.

\begin{lemma}
\label{L_Uniform} If $\mathbf{M\in }\mathcal{M}_{1}$ and $\mathbf{F}(\mathbf{s})\neq \mathbf{Ms},$  then $F_{i}(n;\mathbf{s}%
)\rightarrow 1$ as $n\rightarrow\infty$ uniformly in $i\in\mathbb{N}$ and $%
\mathbf{s}\in(0,1]^{\mathbb{N}}$.
\end{lemma}

\textbf{Proof}. Clearly,
\[
F_{i}(n;\mathbf{s})=\mathbb{E}\left[ \prod_{j\in \mathbb{N}}s_{j}^{Z_{ij}(n)}%
\right] \geq \mathbb{P}\big(\big\|\mathbf{Z}_{i}(n)\big\|_{1}=0\big).
\]%
Recalling Lemma \ref{L_Dichotom} we see that, for each fixed $i\in \mathbb{N}
$
\begin{equation}
\sup_{s\in \lbrack 0,1]^{\mathbb{N}}}\left( 1-F_{i}(n;\mathbf{s})\right)
=\sup_{s\in \lbrack 0,1]^{\mathbb{N}}}Q_{i}(n;\mathbf{s})\leq Q_{i}(n;%
\mathbf{0})\rightarrow 0  \label{hvo5a}
\end{equation}%
as $n\rightarrow \infty $. We now show that convergence in (\ref{hvo5a}) is
uniform over $i\in \mathbb{N}$.

Since
\[
\mathbf{1}-\mathbf{F}(n;\mathbf{s})=\mathbf{1}-\mathbf{F}(1;\mathbf{F}(n-1;%
\mathbf{s}))\leq \mathbf{M}\left( \mathbf{1}-\mathbf{F}(n-1;\mathbf{s}%
)\right) ,
\]%
it follows that for all $N\in \mathbb{N}$
\begin{equation}
Q_{i}(n;\mathbf{s})\leq \sum_{j\in \mathbb{N}}M_{ij}Q_{j}(n-1;\mathbf{s}%
)\leq \sum_{j\leq N}M_{ij}Q_{j}(n-1;\mathbf{s})+\sum_{j>N}M_{ij}.
\label{EstR_i}
\end{equation}

In view of Condition $\mathrm{(iii)}$ describing properties of matrices
belonging to class $\mathcal{M}_{1}$, for any $\varepsilon >0$ there exists
an $N=N(\varepsilon )$ meeting estimate (\ref{le2}). On the other hand, for
any fixed $N$ and $i\in \mathbb{N}$
\begin{equation}
\sum_{j\leq N}M_{ij}\leq \frac{\textstyle{\sum\limits_{j\leq N}M_{ij}u_{j}}}{%
\textstyle{\min\limits_{k\leq N}u_{k}}}\leq \frac{\textstyle{%
\sum\limits_{j\in \mathbb{N}}M_{ij}u_{j}}}{\textstyle{\min\limits_{k\leq
N}u_{k}}}=\frac{u_{i}}{{\min\limits_{k\leq N}u_{k}}}\leq \frac{U}{{%
\min\limits_{k\leq N}u_{k}}}.  \label{ed6}
\end{equation}

Lemma \ref{L_Dichotom} and estimates (\ref{hvo5a})-(\ref{ed6}) imply
\[
\sup_{\mathbf{s}\in \lbrack 0,1]^{\mathbb{N}},\,i\in \mathbb{N}}Q_{i}(n;%
\mathbf{s})\leq \frac{U}{{\min\limits_{k\leq N}u_{k}}}\sup_{j\leq
N}Q_{j}(n-1;\mathbf{0})+\varepsilon \leq 2\varepsilon
\]%
for all sufficiently large $n$.

Lemma \ref{L_Uniform} is proved.

\section{Ratio limit theorem \label{SecRatio}}

The next important theorem is an infinite-dimensional analog of Theorem 1 in
\cite[Ch. VI, \S 1]{Se71}.

Denote $\mathbf{S}:=\left\{ \mathbf{s}\in \lbrack 0,1]^{\mathbb{N}},\mathbf{s%
} \neq \mathbf{1}\right\} $.

\begin{theorem}
\label{t01_splitted} Let $\left\{ \mathbf{Z}_{i}(n),i\in \mathbb{N}\right\}$
be a critical GWBP/$\infty$ with mean matrix $\mathbf{M}\in \mathcal{M}%
_{1}^{0}$ and $\mathbf{F}(\mathbf{s})\neq \mathbf{Ms}.$ Then
\begin{equation}  \label{ted1}
\lim_{n\rightarrow \infty }\sup_{s\in \mathbf{S},\,i\in \mathbb{N}%
}\left\vert \frac{Q_{i}(n;\mathbf{s})}{u_{i}q(n;\mathbf{s})}-1\right\vert =0.
\end{equation}
\end{theorem}

To justify (\ref{ted1}) for the critical GWBP/$\infty $'s we need to
attract, along with the standard hypotheses $\mathrm{(i)}$ for the mean
matrix, additional Conditions $\mathrm{(iii)}$ and $\mathrm{(iv)}$ which
provide the desired uniform convergence in $i\in \mathbb{N}$. These
additional conditions automatically fulfil for the Galton-Watson processes
with finite number of types.

We split the proof of Theorem \ref{t01_splitted} into several lemmas.

For an infinite-dimensional vector $\mathbf{s}=\big(s_{j}\big)_{j\in \mathbb{%
N}}\in \lbrack 0,1]^{\mathbb{N}}$ introduce the notation $\mathbf{s}%
_{j}:=(s_{i}\delta _{ij}+1-\delta _{ij})_{i\in \mathbb{N}}$. Set also
\begin{eqnarray*}
&&F_{ij}(s_{j}):=F_{i}(\mathbf{s}_{j})=\mathbb{E}s_{j}^{Z_{ij}},\ \
Q_{ij}(s_{j}):=Q_{i}(\mathbf{s}_{j})=1-\mathbb{E}s_{j}^{Z_{ij}}, \\
&&\mathcal{N}_{ij}(\mathbf{s}):=\sum_{k=0}^{Z_{ij}-1}s_{j}^{k}\left(
1-\prod_{l=j+1}^{\infty }s_{l}^{Z_{il}}\right) =\frac{1-s_{j}^{Z_{ij}}}{%
1-s_{j}}\left( 1-\prod_{l=j+1}^{\infty }s_{l}^{Z_{il}}\right) , \\
&&N_{ij}(\mathbf{s}):=\mathbb{E}\mathcal{N}_{ij}(\mathbf{s}).
\end{eqnarray*}

\begin{lemma}
\label{loc1} If all elements of the mean matrix $\mathbf{M}=\left(%
\mathbb{E}Z_{ij} \right)_{i,j\in \mathbb{N}}$ are finite then, for each $%
i\in \mathbb{N}$ the following representation is valid
\begin{eqnarray}  \label{res1}
Q_{i}(\mathbf{s}) &=&\sum\limits_{j\in \mathbb{N}}\mathbb{E}\big( %
1-s_{j}^{Z_{ij}}\big)-\sum\limits_{j\in \mathbb{N}}(1-s_{j})\mathbb{E}\left[
\sum\limits_{k=0}^{Z_{ij}-1}s_{j}^{k}\Big(1-\prod\limits_{l=j+1}^{\infty}
s_{l}^{Z_{il}}\Big)\right]  \notag \\
&=&\sum\limits_{j\in \mathbb{N}}Q_{ij}(s_{j})-\sum\limits_{j\in \mathbb{N}%
}(1-s_{j})N_{ij} (\mathbf{s}).
\end{eqnarray}
\end{lemma}

\textbf{Proof.} Using definitions (\ref{dfF1}) we perform a chain of evident
transformations
\begin{eqnarray}  \label{n7}
&&1-\prod_{l=1}^{\infty }s_{l}^{Z_{il}}=\big(1-s_{1}^{Z_{i1}}\big)-\big( %
1-s_{1}^{Z_{i1}}\big)\left( 1-\prod_{l=2}^{\infty }s_{l}^{Z_{il}}\right)
+1-\prod_{l=2}^{\infty }s_{l}^{Z_{il}}  \notag \\
&=&\big(1-s_{1}^{Z_{i1}}\big)-(1-s_{1})\sum_{k=0}^{Z_{i1}-1}s_{1}^{k}
\left(1-\prod_{l=2}^{\infty}s_{l}^{Z_{il}}\right) +\left[ 1-\prod_{l=2}^{%
\infty }s_{l}^{Z_{il}}\right]  \notag \\
&=&\big(1-s_{1}^{Z_{i1}}\big)-(1-s_{1})\mathcal{N}_{i1}(\mathbf{s})+ \left[%
1-\prod_{l=2}^{\infty }s_{l}^{Z_{il}}\right].
\end{eqnarray}

Since all the summands in the chain of identities have finite means, it
follows that
\begin{eqnarray}  \label{n4}
Q_{i}(\mathbf{s}) &=&1-F_{i}(\mathbf{s})=\mathbb{E}\left(
1-\prod_{l=1}^{\infty }s_{l}^{Z_{il}}\right) =\mathbb{E}\left(
1-\prod_{l=2}^{\infty }s_{l}^{Z_{il}}\right)  \notag \\
&+&Q_{i1}(s_{1})-(1-s_{1})N_{i1}(\mathbf{s}).
\end{eqnarray}

Repeating the chain of transformations (\ref{n7}) for the first summand at
the right-hand side of (\ref{n4}) and doing the same for $%
1-\prod_{l=j+1}^{\infty }s_{l}^{Z_{il}}$, where the parameter $j\in \mathbb{N%
}\smallsetminus \{1\}$ sequentially takes values increasing by 1, we get the
desired identity
\[
Q_{i}(\mathbf{s})=\sum\limits_{j\in \mathbb{N}}Q_{ij}(s_{j})-\sum\limits_{j%
\in \mathbb{N}}N_{ij}(\mathbf{s})(1-s_{j}).
\]

The lemma is proved.

\begin{theorem}
\label{t01} Let $\left\{\mathbf{Z}_{i}(n),i\in \mathbb{N}\right\}$ be a
critical GWBP/$\infty$ with mean matrix $\mathbf{M}\in \mathcal{M}_{1}$ and $\mathbf{F}(\mathbf{s})\neq \mathbf{Ms}.$ 
Then
\begin{equation}  \label{ted2}
\lim_{z\uparrow 1}\sup_{i\in \mathbb{N}}\frac{(1-z)\mathbb{E}Z_{i}-\mathbb{E}
(1-z^{Z_{i}})}{(1-z)\mathbb{E}Z_{i}}=0
\end{equation}
and
\[
\lim_{z\uparrow 1}\sup_{i\in \mathbb{N}}\frac{(1-z)\mathbb{E}Z_{i}-
\sum\limits_{j\in \mathbb{N}}\mathbb{E}(1-z^{Z_{ij}})}{(1-z)\mathbb{E}Z_{i}}%
=0.
\]
\end{theorem}

\textbf{Remark.} The statement of the theorem is always true for the
Galton-Watson branching processes with finite number of types. However, this
is not always the case for GWBP/$\infty $'s. It is for this reason we are
forced to include Conditions $\mathrm{(iii)}$ in Definition \ref{d2}.

\textbf{Proof of Theorem \ref{t01}}. Since
\[
\mathbb{E}(1-z^{Z_{i}})\leq \sum\limits_{j\in \mathbb{N}}\mathbb{E}%
(1-z^{Z_{ij}}) \leq \sum\limits_{j\in \mathbb{N}}(1-z)\mathbb{E}Z_{ij}=(1-z)%
\mathbb{E}Z_{i},
\]
it sufficient to prove (\ref{ted2}) only. According to the second part of
Conditions $\mathrm{(iii)}$, for any $\varepsilon >0$ there exists $%
K=K(\varepsilon)$ such that
\begin{equation}  \label{ted7}
0\leq \sup_{i\in\mathbb{N}}M_{i}^{-1}\mathbb{E}\left[
(1-z)Z_{i}-(1-z^{Z_{i}})\,;\,Z_{i}>K\right] \leq \varepsilon \big(1-z\big)/2.
\end{equation}

On the other hand, for any fixed $k\leq K$ the representation
\begin{eqnarray*}
k(1-z)-(1-z^{k}) &=&(1-z)\left( k-\sum_{j=0}^{k-1}z^{j}\right) \\
&\leq &(1-z)\left( K-\sum_{j=0}^{K-1}z^{j}\right)\underset{z\uparrow 1}{=}o%
\big(1-z\big)
\end{eqnarray*}
is valid. Since $\mathbb{P}\left( 0<Z_{i}\leq K\right)\leq \mathbb{E}%
Z_{i}=M_{i}$ for each fixed $K\in \mathbb{N}$, it follows that
\[
0\leq \sup_{i\in \mathbb{N}}M_{i}^{-1}\mathbb{E}\left[%
(1-z)Z_{i}-(1-z^{Z_{i}})\,;\,Z_{i}\leq K\right] \underset{z\uparrow 1}{=}o%
\big(1-z\big).
\]

Thus, for $\varepsilon >0$ and $K=K(\varepsilon)$ selected above one can
find a $\Delta >0$ such that for $0\leq 1-z<\Delta $
\begin{equation}  \label{ted9}
0\leq \sup_{i\in \mathbb{N}}M_{i}^{-1}\mathbb{E}\left[
(1-z)Z_{i}-(1-z^{Z_{i}})\,;\,Z_{i}\leq K\right] \leq \varepsilon \big(1-z%
\big)/2.
\end{equation}

Combining (\ref{ted7}) and (\ref{ted9}) gives (\ref{ted2}).

Theorem \ref{t01} is proved.

Introduce the notation $\mathcal{Q}_{n}=\mathcal{Q}(n;\mathbf{s}):=
\sup\limits_{i\in \mathbb{N}}Q_{i}(n;\mathbf{s})$, $\mathcal{F}_{n}=\mathcal{%
F}(n; \mathbf{s}):=\inf\limits_{i\in \mathbb{N}}F_{i}(n;\mathbf{s})$.

\begin{lemma}
\label{loc11} Let $\left\{ \mathbf{Z}_{i}(n),i\in \mathbb{N}\right\} $ be a
critical GWBP/$\infty $ with mean matrix $\mathbf{M}\in \mathcal{M}_{1}^{0}$ and $\mathbf{F}(\mathbf{s})\neq \mathbf{Ms}.$ Then, for each $i$
\begin{equation}
\sum_{j\in \mathbb{N}}Q_{j}(n;\mathbf{s})N_{ij}\big(\mathbf{F}(n;\mathbf{s})%
\big)=:M_{i}\epsilon _{1,i}(\mathbf{Q}(n;\mathbf{s})),  \label{N21}
\end{equation}%
where
\[
\lim_{n\rightarrow \infty }\sup_{s\in \mathbf{S},\,i\in \mathbb{N}}\frac{%
\left\vert \epsilon _{1,i}(\mathbf{Q}(n;\mathbf{s}))\right\vert }{\mathcal{Q}%
(n;\mathbf{s})}=0.
\]
\end{lemma}

\textbf{Proof.} Recall that $\sup_{i\in \mathbb{N}}M_{i}<\infty $ by (\ref%
{ed04}). By definition
\begin{eqnarray*}
&&Q_{j}(n;\mathbf{s})N_{ij}\big(\mathbf{F}(n;\mathbf{s})\big)=\mathbb{E}%
\left[ \left( 1-F_{j}^{Z_{ij}}(n;\mathbf{s})\right) \left(
1-\prod\limits_{l=j+1}^{\infty }F_{l}^{Z_{il}}(n;\mathbf{s})\right) \right]
\\
&\leq &\mathbb{E}\left[ \big(1-\mathcal{F}_{n}^{Z_{ij}}\big)\left(
1-\prod\limits_{l=j+1}^{\infty }\mathcal{F}_{n}^{Z_{il}}\right) \right] =%
\mathbb{E}\left[ \big(1-\mathcal{F}_{n}^{Z_{ij}}\big)\left( 1-\mathcal{F}%
_{n}^{\sum_{l>j}Z_{il}}\right) \right] .
\end{eqnarray*}

Using the equality $(1-x)(1-y)=1-x+1-y-1+xy$ we see that
\begin{eqnarray}  \label{N24}
&&Q_{j}(n;\mathbf{s})N_{ij}\big(\mathbf{F}(n;\mathbf{s})\big)\leq \mathbb{E} %
\left[ \big(1-\mathcal{F}_{n}^{Z_{ij}}\big)\left( 1-\mathcal{F}%
_{n}^{\sum_{l>j}Z_{il}}\right) \right]  \notag \\
&=&\mathbb{E}\big(1-\mathcal{F}_{n}^{Z_{ij}}\big)+\mathbb{E}\left( 1-
\mathcal{F}_{n}^{\sum_{l>j}Z_{il}}\right)-\mathbb{E}\left(1-\mathcal{F}%
_{n}^{\sum_{l\geq j} Z_{il}}\right) .
\end{eqnarray}

By (\ref{N24}) we conclude that
\begin{eqnarray}
&&0\leq \sum_{j\in \mathbb{N}}Q_{j}(n;\mathbf{s})N_{ij}\big(\mathbf{F}(n;%
\mathbf{s})\big)  \notag \\
&\leq &\sum_{j\in \mathbb{N}}\left( \mathbb{E}\left( 1-\mathcal{F}%
_{n}^{Z_{ij}}\right) +\mathbb{E}\left( 1-\mathcal{F}_{n}^{\sum_{l>j}Z_{il}}%
\right) -\mathbb{E}\left( 1-\mathcal{F}_{n}^{\sum_{l\geq j}Z_{il}}\right)
\right)  \notag \\
&=&\sum_{j\in \mathbb{N}}\mathbb{E}\big(1-\mathcal{F}_{n}^{Z_{ij}}\big)-%
\mathbb{E}\left( 1-\mathcal{F}_{n}^{\sum_{l\in \mathbb{N}}Z_{il}}\right)
\leq \sum_{j\in \mathbb{N}}\mathbb{E}Z_{ij}\mathcal{Q}_{n}-\mathbb{E}\big(1-%
\mathcal{F}_{n}^{Z_{i}}\big)  \notag \\
&=&\mathbb{E}Z_{i}\mathcal{Q}_{n}-\mathbb{E}\big(1-\mathcal{F}_{n}^{Z_{i}}%
\big).  \label{N25}
\end{eqnarray}%
Observe that the series in the second line of (\ref{N25}) uniformly
converges in $i\in \mathbb{N}$ in view of the estimates
\[
\mathbb{E}\left( 1-w^{\sum_{l\geq j}Z_{il}}\right) -\mathbb{E}\left(
1-w^{\sum_{l>j}Z_{il}}\right) \leq \mathbb{E}\left( 1-w^{Z_{ij}}\right)
\]%
and $\mathbb{E}\big(1-w^{Z_{ij}}\big)\leq \mathbb{E}Z_{ij}(1-w)$.

Using (\ref{ted2}) for $x=\mathcal{Q}(n;\mathbf{s})$ we deduce that
\begin{eqnarray*}
0 &\leq &\lim_{n\rightarrow \infty }\sup_{s\in \mathbf{S},\,i\in \mathbb{N}}
\frac{\sum_{j\in \mathbb{N}}Q_{j}(n;\mathbf{s})N_{ij}\big(\mathbf{F}(n;
\mathbf{s})\big)}{M_{i}\mathcal{Q}(n;\mathbf{s})} \\
&\leq &\lim_{n\rightarrow \infty }\sup_{s\in \mathbf{S},\,i\in \mathbb{N}}
\frac{\mathbb{E}Z_{i}\mathcal{Q}(n;\mathbf{s})-\mathbb{E}\left( 1-\mathcal{F}
^{Z_{i}}(n;\mathbf{s})\right) }{M_{i}\mathcal{Q}(n;\mathbf{s})}=0,
\end{eqnarray*}
as desired.

Lemma \ref{loc11} is proved.

\begin{lemma}
\label{oc3} Let $\left\{ \mathbf{Z}_{i}(n),i\in \mathbb{N}\right\} $ be a
critical GWBP/$\infty $ with mean matrix $\mathbf{M}\in \mathcal{M}_{1}$ and $\mathbf{F}(\mathbf{s})\neq \mathbf{Ms}.$ 
Then, for all $i\in \mathbb{N}$,
\[
\sum_{j\in \mathbb{N}}M_{ij}(1-s_{j})-\sum\limits_{j\in \mathbb{N}%
}Q_{ij}(s_{j})=:M_{i}\epsilon _{2,i}(\mathbf{1}-\mathbf{s}),
\]%
where
\[
\lim_{s\in \mathbf{S,}\Vert \mathbf{1}-\mathbf{s}\Vert _{\infty }\rightarrow
0}\sup_{i\in \mathbb{N}}\frac{\left\vert \epsilon _{2,i}(\mathbf{1}-\mathbf{s%
})\right\vert }{\Vert \mathbf{1}-\mathbf{s}\Vert _{\infty }}=0.
\]%
Besides,
\begin{equation}
\sum_{j\in \mathbb{N}}M_{ij}Q_{j}(n-1;\mathbf{s})-\sum\limits_{j\in \mathbb{N%
}}Q_{ij}(F_{j}(n-1;\mathbf{s}))=M_{i}\epsilon _{2,i}(\mathbf{Q}(n-1;\mathbf{s%
})),  \label{F2}
\end{equation}%
where
\[
\lim_{n\rightarrow \infty }\sup_{s\in \mathbf{S},\,i\in \mathbb{N}}\frac{%
\left\vert \epsilon _{2,i}(\mathbf{Q}(n-1;\mathbf{s}))\right\vert }{\mathcal{%
Q}(n-1;\mathbf{s})}=0.
\]
\end{lemma}

\textbf{Proof.} Consider the function
\[
f_{ij}(s_{j}):=M_{ij}(1-s_{j})-Q_{ij}(s_{j})=\mathbb{E}\left[
Z_{ij}(1-s_{j})-1+s_{j}^{Z_{ij}}\right]
\]%
with partial derivative
\[
\frac{\partial f_{ij}(s_{j})}{\partial s_{j}}=\mathbb{E}\left[
Z_{ij}\left(s_{j}^{Z_{ij}-1}-1\right)\right] \leq 0,\ j\in \mathbb{N}.
\]

Two last relations mean that the function $f_{ij}(s_{j})$ is nonincreasing, $%
f_{ij}(1)=0$ and, therefore, $f_{ij}(s_{j})\geq 0$ for $s_{j}\in \lbrack
0,1] $. As a result,
\[
f_{ij}(s_{j})\leq f_{ij}(1-\Vert \mathbf{1}-\mathbf{s}\Vert _{\infty })=
\mathbb{E}\left[ Z_{ij}\Vert \mathbf{1}-\mathbf{s}\Vert _{\infty }-1+\left(
1-\Vert \mathbf{1}-\mathbf{s}\Vert _{\infty }\right)^{Z_{ij}}\right] .
\]
Now the assertion of Lemma \ref{oc3} is an evident corollary of Theorem \ref%
{t01}.

We now deduce estimates for $Q_{i}(n;\mathbf{s})$, $i\in \mathbb{N}$, in
terms of the function $Q_{1}(n;\mathbf{s})$.

\begin{lemma}
\label{L_RR} Let $\left\{ \mathbf{Z}_{i}(n),i\in \mathbb{N}\right\} $ be a
critical GWBP/$\infty$ with mean matrix $\mathbf{M}\in \mathcal{M}_{1}^{0}$ and $\mathbf{F}(\mathbf{s})\neq \mathbf{Ms}.$
Then there exist constants $C_{1}$ and $C_{2}$ such that
\[
Q_{i}(n;\mathbf{s})\leq C_{1}u_{i}Q_{1}(n;\mathbf{s}),\,\forall \,i\in
\mathbb{N},
\]
and
\[
\mathcal{Q}(n;\mathbf{s})\leq C_{2}Q_{1}(n;\mathbf{s})
\]
for all sufficiently large $n$.
\end{lemma}

\textbf{Proof.} Taking into account (\ref{res1}) with $\mathbf{s}$ replaced
by $\mathbf{F}(n-1;\mathbf{s})$, (\ref{N21}) with $\mathbf{F}(n;\mathbf{s})$
replaced by $\mathbf{F}(n-1;\mathbf{s})$, and attracting (\ref%
{F2}) we get
\begin{equation}
Q_{i}(n;\mathbf{s})=\sum_{j\in \mathbb{N}}M_{ij}Q_{j}(n-1;\mathbf{s}%
)+M_{i}\epsilon _{i}(\mathbf{Q}(n-1;\mathbf{s})),  \label{n111}
\end{equation}%
where $\epsilon _{i}(\mathbf{Q}(n-1;\mathbf{s}))=\epsilon _{1,i}(\mathbf{Q}%
(n-1;\mathbf{s}))+\epsilon _{2,i}(\mathbf{Q}(n-1;\mathbf{s}))$ and
\begin{equation}
\lim_{n\rightarrow \infty }\sup_{s\in \mathbf{S},\,i\in \mathbb{N}}\frac{%
\left\vert \epsilon _{i}(\mathbf{Q}(n-1;\mathbf{s}))\right\vert }{\mathcal{Q}%
(n-1;\mathbf{s})}=0.  \label{NeglE_i}
\end{equation}

Multiplying the left and right hand sides of (\ref{n111}) by $v_{i}$ and
summing over $i$ we obtain
\begin{equation}
q(n;\mathbf{s})=\mathbf{vQ}^{T}(n;\mathbf{s})=q(n-1;\mathbf{s})+\epsilon (%
\mathbf{Q}(n-1;\mathbf{s})),  \label{v1}
\end{equation}%
where
\[
\epsilon (\mathbf{Q}(n-1;\mathbf{s})):=\sum_{i\in \mathbb{N}%
}v_{i}M_{i}\epsilon _{i}(\mathbf{Q}(n-1;\mathbf{s})).
\]

In view of (\ref{NeglE_i})
\begin{equation}
\lim_{n\rightarrow \infty }\sup_{s\in \mathbf{S}}\frac{\left\vert \epsilon (%
\mathbf{Q}(n-1;\mathbf{s}))\right\vert }{\mathcal{Q}(n-1;\mathbf{s})}=0.
\label{v1v}
\end{equation}

Representation (\ref{n111}) and estimate (\ref{ed04}) imply, for
sufficiently large $n$ the inequalities
\begin{equation}
Q_{i}(n;\mathbf{s})\leq 2M_{i}\mathcal{Q}_{n-1},\ \ \mathcal{Q}_{n}\leq 2%
\mathfrak{m}\mathcal{Q}_{n-1}.  \label{n112}
\end{equation}

Combining (\ref{n111}) with $n$ sequentially replaced by $n+1,\ldots,$ $n+l$ we
obtain
\begin{equation}
Q_{i}(n+l;\mathbf{s})=M_{i}\epsilon _{i}^{\ast }(n,l;\mathbf{s})+\sum_{j\in
\mathbb{N}}M_{ij}^{(l)}Q_{j}(n;\mathbf{s}),  \label{n114}
\end{equation}%
where
\begin{eqnarray*}
&&M_{i}\epsilon _{i}^{\ast }(n,1;\mathbf{s})=M_{i}\epsilon _{i}(\mathbf{Q}(n;%
\mathbf{s})), \\
&&M_{i}\epsilon _{i}^{\ast }(n,k;\mathbf{s})=\sum\limits_{j\in \mathbb{N}%
}M_{ij}M_{j}\epsilon _{j}^{\ast }(n,k-1;\mathbf{s})+M_{i}\epsilon _{i}(%
\mathbf{Q}(n+k-1;\mathbf{s})),\,k>1.
\end{eqnarray*}%
Clearly,
\[
\lim_{n\rightarrow \infty }\sup_{s\in \mathbf{S},\,i\in \mathbb{N}}\frac{%
\left\vert \epsilon _{i}^{\ast }(n,l;\mathbf{s})\right\vert }{\mathcal{Q}(n;%
\mathbf{s})}=0
\]%
in view of (\ref{NeglE_i}), (\ref{n112}) and (\ref{ed04}).

Since $\mathbf{M}\in\mathcal{M}_{1}^{0}$ by our assumptions, applying the
first inequality in (\ref{ed03}) gives
\begin{equation}  \label{ed004}
M_{ij}^{(n)}=\sum_{k\in \mathbb{N}}M_{ik}M_{kj}^{(n-1)}\leq Cu_{i}\sum_{k\in
\mathbb{N}}v_{k}M_{kj}^{(n-1)}=Cu_{i}v_{j},\ \forall i,j\in \mathbb{N}.
\end{equation}

Further, using the second inequality in (\ref{ed03}) for the respective $m\in%
\mathbb{N}$ we obtain
\begin{equation}  \label{ed005}
M_{1j}^{(n+m)}=\sum_{k\in \mathbb{N}}M_{1k}^{(m)}M_{kj}^{(n)}\geq
c\sum_{k\in \mathbb{N}}v_{k}M_{kj}^{(n)}=cv_{j},\ \forall j\in \mathbb{N}.
\end{equation}

We fix $l\geq m$. Relations (\ref{ed004}) and (\ref{ed005}) imply the
inequalities
\[
M_{ij}^{(l)}\leq Cu_{i}v_{j}\leq \frac{CM_{1j}^{(l)}u_{i}}{c},\ \forall
i,j\in \mathbb{N}.
\]

Hence, using (\ref{n114}) we deduce the estimate
\begin{equation}  \label{ed080}
Q_{i}(n+l;\mathbf{s})\leq \frac{2Cu_{i}}{c}\sum_{j\in \mathbb{N}%
}M_{1j}^{(l)} Q_{j}(n;\mathbf{s}),\ \forall i\in \mathbb{N},
\end{equation}
which, on account of $\left\Vert\mathbf{u}\right\Vert_{\infty }=U<\infty$,
leads to the inequality
\begin{equation}  \label{ed008}
\mathcal{Q}(n+l;\mathbf{s})=\mathcal{Q}_{n+l}\leq \frac{2CU}{c}\sum_{j\in%
\mathbb{N}} M_{1j}^{(l)}Q_{j}(n;\mathbf{s}).
\end{equation}

Finally, using (\ref{n114}) with $i=1$ we may transform (\ref{ed080}) and (%
\ref{ed008}) as
\begin{equation}  \label{ed007}
Q_{i}(n+l;\mathbf{s})\leq \frac{4Cu_{i}}{c}Q_{1}(n+l;\mathbf{s}),\,\forall
\,i\in \mathbb{N},
\end{equation}
and
\begin{equation}  \label{ed009}
\mathcal{Q}_{n+l}=\mathcal{Q}(n+l;\mathbf{s})\leq \frac{4CU}{c}Q_{1}(n+l;%
\mathbf{s}).
\end{equation}

Lemma \ref{L_RR} is proved.

\begin{lemma}
\label{2Sev} If $\mathbf{M}\in \mathcal{M}_{1}^{0}$ and $\mathbf{F}(\mathbf{s})\neq \mathbf{Ms}$ then, for any fixed $n\in
\mathbb{N}$
\begin{equation}  \label{hvo}
\lim_{N\rightarrow \infty }\sup_{s\in \mathbf{S}}\frac{1}{q(n;\mathbf{s})}
\sum_{j=N+1}^{\infty }v_{j}Q_{j}(n;\mathbf{s})=0
\end{equation}
and, for any $\varepsilon>0$ there exist $l_{0}=l_{0}(\varepsilon)\in\mathbb{%
N}$ such that
\begin{eqnarray}
&&\big|Q_{i}(n+l_{0};\mathbf{s})-u_{i}q(n;\mathbf{s})\big|\leq \varepsilon
u_{i}q(n;\mathbf{s}),  \label{nhvo1}\\
&&\big|q(n+l_{0};\mathbf{s})-q(n;\mathbf{s})\big|\leq \varepsilon q(n;%
\mathbf{s})  \label{nhvo2}
\end{eqnarray}
for all $i\in \mathbb{N}$ and $n>n_{0}$.
\end{lemma}

\textbf{Proof.} Clearly,
\[
q(n;\mathbf{s})\geq v_{1}Q_{1}(n;\mathbf{s}),\ \mathbf{s}\in(0,1]^{\mathbb{N}%
}.
\]

On the other hand, estimate (\ref{ed009}) and Conditions (\ref{ed3}) allow
us to deduce for sufficiently large $n$ the estimate
\[
q(n;\mathbf{s})\leq \frac{\textstyle{4CU}}{\textstyle{c}}Q_{1}(n;\mathbf{s}).
\]

Thus, the functions $q(n;\mathbf{s})$ and $Q_{1}(n;\mathbf{s})$ have the
same order as $n\rightarrow\infty$. This fact, estimate (\ref{ed007}) and
convergence of the series $\sum_{i=1}^{\infty}v_{i}$ justify (\ref{hvo}).

Besides,
\[
0<\liminf_{n\rightarrow \infty }\inf_{s\in\mathbf{S}}\frac{q(n;\mathbf{s})} {%
Q_{1}(n;\mathbf{s})}\leq\limsup_{n\rightarrow \infty }\sup_{s\in\mathbf{S}}
\frac{q(n;\mathbf{s})}{Q_{1}(n;\mathbf{s})}<\infty .
\]

We rewrite this relation as
\begin{equation}  \label{hv4}
q(n;\mathbf{s})\mathop{\asymp}\limits_{n\to\infty}Q_{1}(n;\mathbf{s}).
\end{equation}

The same notation will be used in others similar situations. For instance,
in view of (\ref{ed009}) and the definition of $\mathcal{Q}(n+l;\mathbf{s})$
\begin{equation}  \label{hv004}
\mathcal{Q}_{n+l}=\mathcal{Q}(n+l;\mathbf{s})\mathop{\asymp}\limits_{n\to\infty}
Q_{1}(n+l;\mathbf{s})
\end{equation}
for each fixed $l$.

Taking into account Condition $\mathrm{(iv)}$ and relation (\ref{hv004}) we
transform (\ref{n114}) to the form
\begin{eqnarray}  \label{hv333}
&&Q_{i}(n+l;\mathbf{s})-u_{i}q(n;\mathbf{s})-M_{i}\epsilon _{i}^{\ast }(n,l;%
\mathbf{s})  \notag \\
&&\,=\sum_{j=1}^{N}\big(M_{ij}^{(l)}-u_{i}v_{j}\big)Q_{j}(n;\mathbf{s})-
\sum_{j=N+1}^{\infty} u_{i}v_{j}Q_{j}(n;\mathbf{s})+\sum_{j=N+1}^{\infty
}M_{ij}^{(l)}Q_{j}(n;\mathbf{s})  \notag \\
&&\,=:I_{1}(i,N,l,n;\mathbf{s})+I_{2}(i,N,n;\mathbf{s})+I_{3}(i,N,l,n;%
\mathbf{s}).
\end{eqnarray}

Estimates (\ref{ed004}), (\ref{ed009}), (\ref{hv4}), (\ref{hv004}) and the
equality $\mathbf{v1}^{T}=1$ allow us to claim that
\begin{equation}  \label{hv204b}
\lim_{N\rightarrow\infty}\sup_{i\in\mathbb{N},\,n\in\mathbb{N},\,s\in
\mathbf{S}}\left\vert\frac{I_{2}(i,N,n;\mathbf{s})}{u_{i}Q_{1}(n;\mathbf{s})}%
\right\vert =0,
\end{equation}
and
\begin{equation}  \label{hv204}
\lim_{N\rightarrow\infty}\sup_{i\in\mathbb{N},\,l\in\mathbb{N},\,n\in\mathbb{%
N},\,s\in\mathbf{S}}\left\vert \frac{I_{3}(i,N,l,n;\mathbf{s})}{u_{i}Q_{1}(n;%
\mathbf{s})}\right\vert =0
\end{equation}
where $Q_{1}(n;\mathbf{s})$ may be replaced by $q(n;\mathbf{s})$.

We now select an $\varepsilon >0$. In view of (\ref{hv204b}) and (\ref{hv204}%
) there exists $N=N(\varepsilon )$ such that
\begin{equation}
\big|I_{2}(i,N(\varepsilon ),n;\mathbf{s})\big|+\big|I_{3}(i,N(\varepsilon
),l,n;\mathbf{s})\big|\leq 0.25\varepsilon u_{i}q(n;\mathbf{s}).
\label{hv104}
\end{equation}

According to (\ref{ed2}), conditions (\ref{ed3}) and estimate (\ref{ed009}),
there exist $C_{1}$ and $l_{0}=l_{0}(\varepsilon )$ such that, for all $%
l\geq l_{0}$
\begin{equation}
\big|I_{1}(i,N(\varepsilon ),l,n;\mathbf{s})\big|{\leq }N(\varepsilon )C_{1}%
\mathcal{Q}(n;\mathbf{s})\sup_{j\leq N(\varepsilon )}\left\vert M_{ij}^{(l)}{%
-}u_{i}v_{j}\right\vert {\leq }0.25\varepsilon u_{i}q(n;\mathbf{s}).
\label{hv304}
\end{equation}

We know by (\ref{hv4}), (\ref{hv004}) and (\ref{n112}) that $Q_{1}(n+l;%
\mathbf{s})\mathop{\asymp}\limits_{n\to\infty}Q_{1}(n;\mathbf{s})$ for
any fixed $l$. Therefore, for $l=l_{0}$ and $n>n_{2}=n_{2}(\varepsilon )$
the third term at the left-hand side of (\ref{hv333}) may be evaluated as
\begin{eqnarray}  \label{nhv055}
\Big|M_{i}\epsilon ^{\ast }(n,l_{0};\mathbf{s})\Big| =u_{i}\bigg|\frac{M_{i}%
}{u_{i}} \epsilon^{\ast}(n,l_{0};\mathbf{s})\bigg|\leq Cu_{i}\Big|%
\epsilon^{\ast }(n,l_{0};\mathbf{s}) \Big|\leq0.5\varepsilon u_{i}q(n;%
\mathbf{s}).
\end{eqnarray}

Combining estimates (\ref{hv104}), (\ref{hv304}) and (\ref{nhv055}) and
using decomposition (\ref{hv333}) we easily obtain (\ref{nhvo1}).

Since $\mathbf{uv}^{T}$ $=1$, relation (\ref{nhvo2}) immediately follows
from (\ref{nhvo1}).

Lemma \ref{2Sev} is proved.

\textbf{Proof of Theorem \ref{t01_splitted}.} Clearly, estimates (\ref{nhvo1}) and (\ref{nhvo2}) imply (\ref{ted1}).

Theorem \ref{t01_splitted} is proved.

\section{Proof of Theorem \protect\ref{t02} \label{SecFin}}

We set
\[
B(n;\mathbf{s}):=\mathbf{vQ}^{T}(n;\mathbf{s})-\mathbf{vQ}^{T}(\mathbf{F}(n;%
\mathbf{s}))
\]%
and first prove an infinite-dimensional analog of Lemma 2 in \cite{V1977}.

Recall that $\Phi (x)=x-\mathbf{v}\mathbf{Q}^{T}(\mathbf{1}-x\mathbf{u})
=x^{1+\alpha }\ell (x)$ for $xU\leq1$ as $x\rightarrow +0$ by (\ref{V2}) and
(\ref{vat1}).

\begin{lemma}
\label{vaL2} If the conditions of Theorem \ref{t02} are valid then
\begin{equation}  \label{vat3}
\lim_{n\rightarrow \infty }\sup_{s\in \mathbf{S}}\left\vert \frac{\textstyle{%
B(n;\mathbf{s})}} {\textstyle{\Phi (q(n;\mathbf{s}))}}-1\right\vert =0.
\end{equation}
\end{lemma}

The proof of (\ref{vat3}) coincides almost literally with the proof of Lemma
2 in \cite{V1977} and we give it here to only keep the integrity of the
presentation.

Introduce the function
\[
B(\mathbf{s}):=\mathbf{v}(\mathbf{1}-\mathbf{s})^{T}-\mathbf{vQ}^{T}(\mathbf{%
s}).
\]

Clearly, for $\mathbf{s}=(s_{1},s_{2},...)\in \mathbf{S}$
\begin{eqnarray*}
\frac{\partial B(\mathbf{s})}{\partial s_{i}} &=&-v_{i}-\sum_{j=1}^{%
\infty}v_{j} \frac{\partial Q_{j}(\mathbf{s})}{\partial s_{i}}=
-v_{i}+\sum_{j=1}^{\infty }v_{j} \frac{\partial F_{j}(\mathbf{s})}{\partial
s_{i}} \\
&\leq &-v_{i}+\sum_{j=1}^{\infty}v_{j}\mathbb{E}Z_{ji}=-v_{i}+v_{i}=0.
\end{eqnarray*}
Thus, $B(\mathbf{s})$ is monotone decreasing with respect to each argument
of $\mathbf{s}$. By Theorem \ref{t01_splitted}, for any $\varepsilon >0$ one
can find $N=N(\varepsilon)$ such that
\[
(1-\varepsilon )u_{i}q(n;\mathbf{s})\leq Q_{i}(n;\mathbf{s})\leq
(1+\varepsilon )u_{i}q(n;\mathbf{s})
\]
for all $n\geq N$ and all $i\in\mathbb{N}$ and $\mathbf{s}\in\mathbf{S}$.
Therefore, for $n\geq N$
\[
B(\mathbf{1}-(1-\varepsilon )q(n;\mathbf{s})\mathbf{u})\leq B(\mathbf{1}-
\mathbf{Q}(n;\mathbf{s}))\leq B(\mathbf{1}-(1+\varepsilon )q(n;\mathbf{s})%
\mathbf{u}).
\]

Since $B(\mathbf{1}-\mathbf{Q}(n;\mathbf{s}))=B(n;\mathbf{s})$ and $B(%
\mathbf{1}- x\mathbf{u})=\Phi \left( x\right) $, it follows that, for $n\geq
N$
\[
\Phi \left((1-\varepsilon)q(n;\mathbf{s})\right)\leq B(n;\mathbf{s})\leq
\Phi \left((1+\varepsilon)q(n;\mathbf{s})\right).
\]

By our conditions, $\ell(x)$ is a slowly varying function as $x\rightarrow
+0 $. Therefore (see, for instance, Theorem 1.1 \cite[Ch. 1, \S 1.2]{Sen76}%
),
\[
\frac{\ell (cx)}{\ell (x)}\rightarrow 1
\]
as $x\rightarrow+0$ uniformly in $c\in\lbrack a,b]$, $0<a<b<\infty$. Fix $%
\varepsilon_{0} \in(0,1)$. Since $q(n;\mathbf{s})\leq q(n;\mathbf{0})$ and $%
\lim_{n\rightarrow\infty}q(n; \mathbf{0})=0$, we see that
\[
\lim_{n\rightarrow\infty}\sup_{\mathbf{s}\in\mathbf{S}}\frac{%
\Phi((1\pm\varepsilon) q(n;\mathbf{s}))}{\Phi(q(n;\mathbf{s}))}{=}%
(1\pm\varepsilon )^{1+\alpha} \lim_{n\rightarrow\infty}\sup_{\mathbf{s}\in%
\mathbf{S}}\frac{\ell((1\pm\varepsilon) q(n;\mathbf{s}))}{\ell(q(n;\mathbf{s}%
))}{=}\left(1\pm\varepsilon\right)^{1+\alpha }.
\]
By letting $\varepsilon \rightarrow +0$ we easily deduce (\ref{vat3}).

The next statement is an infinite-dimensional of Lemma 3 in \cite{V1977}.

\begin{lemma}
\label{vaL3} Let the conditions of Theorem \ref{t02} be valid. Then
\begin{equation}  \label{vat4}
\lim_{n\rightarrow \infty }\sup_{\mathbf{s}\in \mathbf{S}}\left\vert \frac{%
\textstyle{\Phi(q(n+1,\mathbf{s}))}}{\textstyle{\Phi(q(n;\mathbf{s}))}}%
-1\right\vert=0.
\end{equation}
\end{lemma}

\textbf{Proof.} In view of $\Phi (x)=x^{1+\alpha }\ell (x)$, to demonstrate
the validity of (\ref{vat4}) it is sufficient to show that
\[
\lim_{n\rightarrow \infty }\sup_{\mathbf{s}\in \mathbf{S}}\left\vert \frac{%
\textstyle{q(n+1,\mathbf{s})}}{\textstyle{q(n;\mathbf{s})}}-1\right\vert =0.
\]
It remains to observe that the desired estimate is a corollary of (\ref{v1}%
), (\ref{v1v}), (\ref{hv4}) and (\ref{hv004}).

Lemma \ref{vaL3} is proved.

\textbf{Proof of Theorem \ref{t02}}. Using Lemma \ref{vaL2} we write
\[
B(k;\mathbf{s})={q(k;\mathbf{s})-q(k+1;\mathbf{s})=\Phi (q(k;\mathbf{s}))}%
\left( 1+\varepsilon (k;\mathbf{s})\right) ,
\]%
where
\[
\lim_{k\rightarrow \infty }\sup_{\mathbf{s}\in \mathbf{S}}\left\vert
\varepsilon (k;\mathbf{s})\right\vert =0.
\]%
Hence, setting $q(0;\mathbf{s}):=\mathbf{v}\left( \mathbf{1}-\mathbf{s}%
\right) ^{T}$ we obtain that
\[
\sum_{k=0}^{n}\frac{q(k;\mathbf{s})-q(k+1;\mathbf{s})}{\Phi (q(k;\mathbf{s}))%
}=n\left( 1+\varepsilon _{1}(n;\mathbf{s})\right) ,
\]%
where $\lim_{n\rightarrow \infty }\sup_{\mathbf{s}\in \mathbf{S}}\left\vert
\varepsilon _{1}(n;\mathbf{s})\right\vert =0$. Lemma \ref{vaL3} and
monotonicity of $\Phi (x)$ in $x$ allow us to rewrite the previous relation
as
\begin{equation}
\int_{q(n;\mathbf{s})}^{q(0;\mathbf{s})}\frac{dx}{\Phi (x)}=n(1+\varepsilon
_{2}(n;\mathbf{s})),  \label{Integ_Phi}
\end{equation}%
where $\lim_{n\rightarrow \infty }\sup_{\mathbf{s}\in \mathbf{S}}\left\vert
\varepsilon _{2}(n;\mathbf{s})\right\vert =0$. Letting $\mathbf{s}=\mathbf{0}
$ and recalling that $\Phi (x)=x^{1+\alpha }\ell (x)$ as $x\rightarrow +0$
we deduce by the properties of regularly varying functions (see Theorem 1
\cite[Ch. VIII, Section 9]{Fe2}) that
\[
q^{\alpha }(n)\ell (q(n))\sim (\alpha n)^{-1},\ n\rightarrow \infty ,
\]%
or (see property $5^{\circ }$ \cite[Ch. 1, Section 1.5]{Sen76})
\[
q(n)=n^{-1/\alpha }\ell _{1}(n)
\]%
for a function $\ell _{1}(n)$ slowly varying as $n\rightarrow \infty $. This
proves (\ref{ted22}).

Relation (\ref{ted2222}) follows from (\ref{ted22}) and Theorem \ref%
{t01_splitted}.

We now prove (\ref{tv1_new}) and (\ref{ted11}). Let $\boldsymbol{\lambda}
=\left(\lambda_{j}\right)_{j\in\mathbb{N}}$ be an infinite-dimensional
vector with nonnegative bounded components. Set
\begin{equation}  \label{0d2}
s_{i}=s(n;\lambda_{i})=\exp\left\{ -\lambda_{i}q(n)\right\} ,\ i=1,2,...,
\end{equation}
and put $q(0;\mathbf{s}):=\mathbf{v}(\mathbf{1}-\mathbf{s})^{T}$. Using the
relation $\Phi (q(n))/q(n)\sim (\alpha n)^{-1}$, $n\rightarrow\infty$, and
making the change of variables $x\rightarrow zq(n)$ we deduce from (\ref%
{Integ_Phi}) the representation
\begin{equation}  \label{d3}
\int\limits_{q(n;\mathbf{s})/q(n)}^{q(0;\mathbf{s})/q(n)}\frac{\textstyle{%
\Phi (q(n))}d{z}} {\textstyle{\Phi(zq(n))}}=\frac{1+\varepsilon_{2}(n;%
\mathbf{s})}{\alpha},
\end{equation}
where $\varepsilon_{2}(n;\mathbf{s})\rightarrow0$ as $n\rightarrow\infty$.
Note that, as $n\rightarrow\infty$
\begin{equation}  \label{PhiUnif}
\frac{\Phi (q(n))}{\Phi (zq(n))}\rightarrow \frac{1}{z^{1+\alpha }}
\end{equation}
uniformly in $z$ from any finite interval $0<a\leq z\leq b<\infty$.

By assumption, the components of $\boldsymbol{\lambda }$ are bounded and $%
\mathbf{v1}^{T}=1$. Therefore,
\[
\lim_{n\rightarrow \infty }\frac{q(0;\mathbf{s})}{q(n)}=(\mathbf{v},%
\boldsymbol{\lambda }).
\]%
Since the right-hand side of (\ref{d3}) has a limit as $n\rightarrow \infty $%
, the same is true for the left-hand side. Consequently,
\[
\lim_{n\rightarrow \infty }q(n;\mathbf{s})/q(n)=:1-\phi (\boldsymbol{\lambda
})
\]%
also exists. Moreover, this limit is strictly positive. Indeed, if it would
be not the case then the integral at the left-hand side of (\ref{d3}) would
be divergent in view of (\ref{PhiUnif}). Using (\ref{PhiUnif}) once again,
passing to the limit in (\ref{d3}) and performing integration we obtain
\[
(1-\phi (\boldsymbol{\lambda }))^{-\alpha }-(\boldsymbol{v},\boldsymbol{%
\lambda })^{-\alpha }=1
\]%
or
\[
\phi (\mathbf{\lambda })=1-\big(1+(\mathbf{v},\boldsymbol{\lambda }%
)^{-\alpha }\big)^{-1/\alpha }.
\]%
Finally, recalling Theorem \ref{t01_splitted} and selecting the same $%
\mathbf{s}$ as in (\ref{0d2}) we obtain
\begin{eqnarray*}
\lim_{n\rightarrow \infty }\mathbb{E}\left[ e^{-\left( \boldsymbol{\lambda },%
\mathbf{Z}(n)\right) q(n)}\big|\mathbf{Z}(n)\neq \mathbf{0},\mathbf{Z}(0)=%
\mathbf{e}_{i}\right] &=&1-\lim_{n\rightarrow \infty }\frac{Q_{i}(n;\mathbf{s%
})}{Q_{i}(n)} \\
&=&1-\lim_{n\rightarrow \infty }\frac{q(n;\mathbf{s})}{q(n)}=\phi (%
\boldsymbol{\lambda }).
\end{eqnarray*}%
The last is equivalent to (\ref{tv1_new}) which, in turn, implies (\ref%
{ted11}).

Theorem \ref{t02} is proved.

\end{document}